\documentclass{amsart}
\usepackage[all]{xy}
\usepackage{amssymb}
\newcommand{\BQ}{\mathbb Q}
\newcommand{\BC}{\mathbb C}
\newcommand{\BN}{\mathbb N}
\newcommand{\BF}{\mathbb F}
\newcommand{\BZ}{\mathbb Z}

\newcommand{\cH}{\mathcal H}
\newcommand{\cA}{\mathcal A}
\newcommand{\cD}{\mathcal D}
\newcommand{\cB}{\mathcal B}

\newcommand{\cF}{\mathcal F}
\newcommand{\cG}{\mathcal G}
\newcommand{\cC}{\mathcal C}
\newcommand{\cP}{\mathcal P}

\newcommand{\uc}{{\bf c}}
\newcommand{\bc}{\uc}
\newcommand{\bI}{{\bf I}}
\newcommand{\be}{{\bf 1}}
\newcommand{\bo}{{\bf 0}}
\newcommand{\bX}{{\bf X}}
\newcommand{\bY}{{\bf Y}}
\newcommand{\bZ}{{\bf Z}}
\newcommand{\bK}{{\bf K}}
\newcommand{\cM}{\mathcal M}
\newcommand{\fI}{\mathbb I}
\newcommand{\Gr}{\mathfrak Gr}
\newcommand{\fM}{\mathfrak M}
\newcommand{\Fl}{\mathfrak Fl}
\newcommand{\G}{{}^LG}
\newcommand{\hW}{\hat W}
\newcommand{\Res}{\mbox{Res}}
\newcommand{\iHom}{\underline{\mbox{Hom}}}
\newcommand{\sq}{$\square$}
\newcommand{\iso}{\widetilde\longrightarrow}
\newcommand{\imbed}{\hookrightarrow}
\newcommand{\Htil}{\tilde H}
\newcommand{\Gm}{{\mathbb G}_m}
\newcommand{\pr}{{\rm pr }}

\title{On tensor categories attached to cells in affine Weyl groups II}
\author{Roman Bezrukavnikov} \email{roman@math.uchicago.edu}
\author{Viktor Ostrik} \email{ostrik@math.mit.edu}
\date{February 2001}
\thanks{The first author was partially supported by
an NSF grant; part of the work
was done while he was employed by the Clay Institute}
\begin{document}
\maketitle
\section{Introduction}
Let $R$ be a root system. Let $W$ be the corresponding
affine Weyl group, and let $\hat W$
be an extended affine Weyl group. Let $\cH$ (respectively $\hat \cH$) be the
corresponding Hecke algebras. George Lusztig defined an asymptotic version of
the  Hecke
algebra, the ring $J$, see \cite{Lc}. By definition the ring $J$ is a direct
sum $J=\bigoplus_{\uc}J_{\uc}$ where summation is over the set of
{\em two-sided cells} in the affine Weyl group. Further, G.~Lusztig proved that
the set of two-sided cells in $W$ is bijective to the set of unipotent
conjugacy classes in an algebraic
group over $\BC$ with root system $R$, see \cite{Lc} IV. Moreover, he proposed
a Conjecture describing rings $J_{\uc}$ in terms of convolution algebras, see
\cite{Lc} IV, 10.5 (a), (b). This Conjecture was verified in many cases by
Nanhua Xi, see \cite{X1, X2, X3}. In this note we give a
 more conceptual proof of
all previously known results. Our proof also works in some new cases. In
general, we prove a statement (see Theorem 4 below) which is weaker than 
Lusztig's Conjecture.

The proof relies on many results of G.Lusztig in \cite{Lc}. Our new essential
tool is the theory of {\em central sheaves} on affine flag manifold due to
A.~Beilinson, D.~Gaitsgory, R.~Kottwitz, see \cite{G}. One of us used this
theory to prove a part of Lusztig's Conjecture, see \cite{B}.

We would like to thank George Lusztig for useful conversations.


\section{Recollections}
\subsection{Notations} Let $G$ be an algebraic reductive connected group over
the field of $l-$adic numbers $\bar \BQ_l$. Let $X$ denote the weight lattice
of $G$ and let $R\subset X$ denote the root system of $G$. Let $W_f$ denote
the Weyl group of $G$ and let $\hat W$ be the extended Weyl group, that is the
semidirect product of $W_f$ and $X$. Let $l: \hat W \to \BZ$ be the length
function. Let $W\subset \hat W$ be the affine Weyl group, that is subgroup
generated by $W_f$ and $R\subset X$. Let $S=\{ s\in W | l(s)=1\}$ be the set
of simple reflections. It is well known that pair $(W, S)$ is a Coxeter system.

It is well known that any right $W_f-$coset in $\hW$ contains unique shortest
element. Let $\hW^f\subset \hW$ denote the subset of such representatives,
so the set $\hW^f$ is in natural bijection with $\hW/W_f$.

\subsection{Affine Hecke algebra} Let $\cA =\BC [v, v^{-1}]$. The affine Hecke
algebra $\hat \cH$ is a free $\cA-$module with basis $H_w (w\in \hat W)$ with
an associative $\cA-$algebra structure defined by $H_wH_{w'}=H_{ww'}$ if
$l(ww')=l(w)+l(w')$ and $(H_s+v^{-1})(H_s-v)=0$ if $s\in S$. The algebra $\hat
\cH$ is endowed with the
{\em Kazhdan-Lusztig basis} $C_w, w\in \hat W$, see e.g.
\cite{Lc} IV 1.1. Let $h_{x,y,z}\in \cA$ be the structure constants of $\hat
\cH$ with respect to this basis, that is
$$
C_xC_y=\sum_{z\in W}h_{x,y,z}C_z.
$$

We say that (left, right or two-sided) ideal $I\subset \hat \cH$ is KL-ideal
if it admits an $\cA-$basis consisting of some elements $C_x$. For $x,y\in
\hat W$ we write $x\le_L y$ (resp. $x\le_Ry, x\le_{LR}y$) if left (resp.
right, two-sided) KL-ideal generated by $x$ is contained in left (resp. right,
two-sided) KL-ideal generated by $y$, cf. \cite{Licm}. The relations
$\le_L, \le_R, \le_{LR}$ are preorders. Let $\sim_L, \sim_R, \sim_{LR}$ be the
associated equivalence relations. The
corresponding equivalence classes are called
left, right and two-sided cells, see {\em loc. cit.} Each two-sided cell is a
union of left (resp. right) cells. The map $w\mapsto w^{-1}$ induces a
 bijection
of the set of left cells to the set of right cells. This map induces identity
on the set of two-sided cells.

A deep Theorem due to G.~Lusztig (see \cite{Lc} IV 4.8) states that the set of
two-sided cells is bijective to the set of unipotent orbits in $G$.


\subsection{Asymptotic Hecke algebra $J$}
There are well defined functions
$a: \hat W\to \BN, \gamma :\hat W\times \hat W\times \hat W\to \BN$ such that
$$
v^{a(z)}h_{x,y,z}-\gamma_{x,y,z^{-1}}\in v\BZ [v]\; \mbox{for all}\; x,y,z\in
\hat W $$
and such that for any $z\in \hat W$ there exist $x,y\in \hat W$ with
$\gamma_{x,y,z}\ne 0$.
The function $a$ is constant on two-sided cells, see \cite{Lc} I 5.4.

Let $J$ be a free $\BZ-$module with basis $t_x, x\in \hat W$. It has a unique
structure of an associative $\BZ-$algebra such that
$t_xt_y=\sum_{z\in \hW}\gamma_{x,y,z}t_{z^{-1}}$, see \cite{Lc} II. It has a
unit element $\sum_{t\in \cD}t_d$ where the summation is over the set
$\cD\subset W$ of {\em distinguished involutions}, see {\em loc. cit.}
Each left (resp. right) cell contains exactly one element of $\cD$.
For any two-sided cell $\bc$ let $J_{\bc}\subset J$ be the $\BZ-$submodule
generated by $t_x, x\in \bc$. The submodule $J_{\bc}$ is in fact a subalgebra;
moreover $J_{\bc}\cdot J_{\bc'}=0$ for $\bc \ne \bc'$, see \cite{Lc} II,
hence $J=\bigoplus_{\bc} J_{\bc}$.

We will use many times the following characterization of cells due to
G.~Lusztig: $w\sim_L w'$ if and only if $t_wt_{{w'}^{-1}}\ne 0$, see \cite{Ll}
3.1 (k).

Algebras $J_{\bc}$ are examples of {\em based algebras}, that is algebras over
$\BZ$ endowed with a basis over $\BZ$ such that the structure constants in this
basis are nonnegative integers. Another example of a based algebra can be
constructed as follows: let $F$ be a reductive algebraic group acting on the
finite set $X$; then the Grothendieck group $K_F(X\times X)$ of the category of
$F-$equivariant coherent sheaves on $X\times X$ is a based algebra with the
basis given by classes of irreducible $F-$bundles and multiplication given by
convolution, see \cite{Lc} IV 10.2.

Assume for a moment that group $G$ is simply connected. For any two-sided cell
$\bc$ let $u_{\bc}$ be the unipotent element in $G$ corresponding to $\bc$
under Lusztig's bijection \cite{Lc} IV 4.8. Let $F_{\bc}$ be the Levi factor
of the centralizer $Z_G(u_{\bc})$ of $u_{\bc}$ in $G$. In \cite{Lc} IV 10.5
G.~Lusztig conjectured that there exists a finite set $\bX$ endowed with an
action of $F_{\bc}$ such that the based algebras $J_{\bc}$ and $K_{F_{\bc}}
(\bX \times \bX)$ are isomorphic as based algebras, that is the isomorphism
respects bases. The aim of this note is to prove a weak form of this
Conjecture; more precisely, we replace finite $F_{\bc}-$set by a somewhat more
general object --- finite $F_{\bc}-$set of {\em centrally extended points},
see below.

\subsection{} \label{peresechenie}
We will need the following well known

{\bf Lemma.} {\em Let $\Gamma_1$ and $\Gamma_2$ be two left cells lying in
the same two-sided cell. Then the
intersection $\Gamma_1\cap (\Gamma_2)^{-1}$ is non empty.}

{\bf Proof.} Let $w\in \Gamma_1$ and $w'\in \Gamma_2^{-1}$.
By \cite{Ll} 3.1 (l) $w\sim_{LR}w'$ if and only if there exists
$y\in \hW$ such that $t_wt_yt_{w'}\ne 0$. We see from the characterization of
left cells above that $w\sim_Ly^{-1}$ and $y\sim_L{w'}^{-1}$. Thus
$y^{-1}\in \Gamma_1\cap \Gamma_2^{-1}$. \sq

\section{Affine flags}
\subsection{Notations}
Let $\G$ be a split reductive algebraic group over
$\BZ$ which is Langlands dual to $G$. To $\G$ one associates the following
``loop objects'' defined over $\BZ$: the (inifinite type) group schemes
$\bK_\BZ$ of maps from a formal disc to $\G$, and  the Iwahori
group $\bI_\BZ$  of maps whose value at the origin lies
in a fixed Borel; and ind-schemes $\Fl_\BZ$ (the affine flag variety),
and $\Gr_\BZ$ (the affine Grassmanian). For a field ${\mathfrak k}$ we have
$\bK({\mathfrak k})=\G(O)$, $\bI({\mathfrak k})=I$, $\Gr({\mathfrak k})
=\G (F)/\G (O)$ and $\Fl ({\mathfrak k})=
\G (F)/I$ where $F={\mathfrak k}((t))$, $O={\mathfrak k}[[t]]$, and
$I\subset \G (O)$ is an Iwahori subgroup.

We fix a field ${\mathfrak k}$ which is either $\overline{\BF}_p$
or complex numbers; we change scalars from $\BZ$ to ${\mathfrak k}$ (and drop
the subscript $\BZ$). By the (derived) category of sheaves we will
mean either the (derived) category of $l$-adic sheaves, $l\ne char({\mathfrak
k})$, or  the (derived) category of constructible sheaves
on the  complex variety for ${\mathfrak k}=\BC$. We will denote
$\overline {\BQ}_l$ by $\BC$ in the first case.

The orbits of $\bI$ on $\Fl$, $\Gr$ are finite dimensional and isomorphic
to affine spaces; it is well known that orbits (called {\em Schubert cells})
are labelled by elements of $\hat W$ for $\Fl$ and $\hat W/W_f$ for $\Gr$.
For $w\in \hat W$ (respectively $w\in \hat W/W_f$) let $\Fl_w$, $\Gr_w$ be
the corresponding Schubert cells.

Let $D^I$ be the
 $\bI-$equivariant derived category of sheaves on $\Fl$, and
let $\cP^I\subset D^I$ be the full subcategory of perverse sheaves. The
convolution product defines a functor $* : D^I\times D^I\to D^I$;
moreover, $*$ is equipped with a
natural associativity constraint (cf. e.g. \cite{BGK}, \S 1.1.2-1.1.3,
or \cite{BD}, \S 7.6.1, p. 260). 

Let $j_w: \Fl_w\to \Fl$ be the natural inclusion and let $L_w=j_{w!*}
(\underline{\BQ_l}[\dim \Fl_w])$, where $\BQ_l$ is the constant sheaf.
Simple objects in $\cP^I$ are exhausted by $L_w, w\in \hat W$.

{\bf Remark.} Following the standard yoga
one can consider  the ``graded'' versions of $D^I_{mix}$, $\cP^I_{mix}$ of
$D^I$, $\cP^I$; here $D^I_{mix}$, $\cP^I_{mix}$ are
subcategories in the derived category of mixed $l$-adic sheaves
if ${\mathfrak k}$ is of finite characteristic, and they are objects of the
(derived) category of mixed Hodge
$D$-modules if ${\mathfrak k}=\BC$. The convolution on  $D^I_{mix}$ is defined.
It provides $D^I_{mix}$ with the structure of a monoidal category, and thus
equips  the Grothendieck group $K(D^I_{mix})$ with an algebra structure;
this algebra is isomorphic to
$ \cH$.

We will not use this theory below; however, it underlies the relation
between the categories considered in this note and
affine Hecke
algebras. Also, since the set of
representations of an affine
 Hecke algebra injects into the set of representations
of the corresponding
$p$-adic group $\G\left(\BF_q((t)) \right)$, appearance of
the Langlands dual group in the statements below is a manifestation of
the geometric Langlands duality.

Notice also that mixed sheaves are used in \cite{B} (in the proof of
Theorem 2); the results of this note are based on those of \cite{B}.

\subsection{Central sheaves} \label{Gait}
Recall the following definition, see e.g. \cite{Kas}

{\bf Definition.} {\em Let $\cA$ be a monoidal category, and $\cB$ be a
tensor (symmetric monoidal) category. A central functor from $\cB$ to $\cA$
is a monoidal functor $F:\cB \to \cA$ together with a functorial isomorphism
$$
\sigma_{X,Y}: F(X)\otimes Y\cong Y\otimes F(X)
$$
fixed
for all $X\in \cB,\; Y\in \cA$ subject to the following compatibilities:}

(i) {\em For $X, X'\in \cB$ the isomorphism $\sigma_{X, F(X')}$
coincides with the isomorphism induced by commutativity constraint in $\cB$.}

(ii) {\em For $Y_1, Y_2\in \cA$ and $X\in \cB$ the composition
$$
F(X)\otimes Y_1\otimes Y_2\stackrel{\sigma_{X,Y_1}\otimes id}{\longrightarrow}
Y_1\otimes F(X) \otimes Y_2 \stackrel{id \otimes \sigma_{X,Y_2}}
{\longrightarrow} Y_1\otimes Y_2 \otimes F(X)
$$
coincides with $\sigma_{X, Y_1\otimes Y_2}$.}

(iii) {\em For $Y\in \cA$ and $X_1, X_2\in \cB$ the composition
$$
F(X_1\otimes X_2)\otimes Y\cong
F(X_1)\otimes F(X_2)\otimes Y\stackrel{id\otimes \sigma_{X_2,Y}}
{\longrightarrow} F(X_1)\otimes Y \otimes F(X_2)\stackrel{\sigma_{X_1,Y}
\otimes id}{\longrightarrow}
$$
$$
Y\otimes F(X_1)\otimes F(X_2)\cong Y\otimes F(X_1\otimes X_2)
$$
coincides with $\sigma_{X_1\otimes X_2, Y}$.}

Let $\cP_{\Gr}$ be the category of $\bK-$ equivariant perverse sheaves on
$\Gr$. The convolution endows $\cP_{\Gr}$ with monoidal structure and this
structure naturally extends to a structure of a commutative rigid tensor
category with a fiber functor, and this category is equivalent to
$Rep(G)$, see \cite{G, MV}; \cite{BD},
\S 5.3, pp 199--215. We will identify $Rep(G)$ with $\cP_{\Gr}$.

In \cite{G} a functor $Z: Rep(G)=\cP_{\Gr}\to \cP^I(\Fl)$ was constructed.
It 
enjoys the following properties:

(i) We have a natural isomorphism of functors
$\pi_*\circ Z\cong id$, where $\pi : \Fl \to \Gr$ is the projection.

(ii) For $\cF \in \cP_{\Gr}$, $\cG \in \cP^I$ we have $\cG *Z(\cF )\in \cP^I$.

(iii) $Z$ is endowed with a natural structure of
a central functor from the tensor category $\cP_{\Gr}$ to the
monoidal category $D^I$.

(iv) A unipotent automorphism (monodromy) $\fM$ of the tensor functor $Z$ is
given; centrality isomorphism from (iii) commutes with $\fM$.

\subsection{Monoidal category $\cA_{\bc}$}
For a subset $S\subset W$ let $\cP^I_S$ denote the Serre subcategory of
$\cP^I$ with simple objects $L_w,\; w\in S$. Let $\hat W_{\le \bc}=
\bigcup_{\bc' \le_{LR} \bc}\bc'$ and $\hat W_{< \bc}=\bigcup_{\bc' <_{LR}
\bc}\bc'$. We abbreviate $\cP^I_{\le \bc}=\cP^I_{\hat W_{\le \bc}}$ and
$\cP_{<\bc}=\cP_{\hat W_{<\bc}}$. Let $\cP^I_{\bc}$ denote the Serre
quotient category $\cP^I_{\le \bc}/\cP^I_{<\bc}$.

For any object $X\in D^I$ and integer $i$ let $H^i(X)\in \cP^I$ denote
$i-$th perverse cohomology. For any $X, Y\in \cP^I_{\bc}$ let us define
truncated convolution $X\bullet Y\in \cP^I_{\bc}$ by $X\bullet Y=H^{a(\bc )}
(X*Y) \mod \cP^I_{<\bc}$. Let $\cM_{\bc}$ be the full subcategory of
$\cP^I_{\bc}$ consisting of semisimple objects. It follows from the
Decomposition Theorem \cite{BBD} that the functor $\bullet$ preserves
category $\cM_{\bc}$. The fact the convolution
of pure perverse sheaves is pure implies (see \cite{Lt}, 2.6)
that the Grothendieck group $K(\cM_{\bc})$
with the multiplication induced by $\bullet$ is isomorphic to the algebra
$J_{\bc}$. In \cite{Lt} a natural associativity constraint was
constructed for $\bullet$. Let $\fI_{\bc}=\bigoplus_{d\in \bc \cap \cD}L_d\in
\cM_{\bc}$ (recall that $\cD$ is the set of distinguished involutions). It is
clear that $\fI_{\bc}\bullet X\simeq X\bullet \fI_{\bc}\simeq X$
for any $X\in \cM_{\bc}$. Thus a choice of an isomorphism
$\fI_{\bc}\bullet \fI_{\bc}\to \fI_{\bc}$ defines a structure of a monoidal
category on $\cM_{\bc}$, see \cite{Lt}. We will fix such a choice for the rest
of this paper.

Let $\cA_{\bc}$ be the full subcategory of $\cP^I_{\bc}$ consisting of all
subquotients of $L_w*Z(\cF ) \mod \cP^I_{<\bc}$ where $w\in \bc$ and
$\cF \in \cP_{\Gr}$. The following Proposition is proved in \cite{B},
Proposition 2.

{\bf Proposition.} {\em Restriction of $\bullet$ to $\cA_{\bc}$ takes
values in $\cA_{\bc}$, is exact in each variable, and it equips $\cA_{\bc}$
with a structure of a monoidal category with unit object $\fI_{\bc}$.}

It is clear from the definitions that Lusztig's category $\cM_{\bc}$ is a
monoidal subcategory of $\cA_{\bc}$ consisting of semisimple objects in
$\cA_{\bc}$.

\subsection{Some results from \cite{B}} \label{Bezr}
Let $d\in \bc$ be a Duflo involution.
Let $\cA_d\subset \cA_{\bc}$ be the full subcategory consisting of all
subquotients of $L_d*Z(\cF ),\; \cF \in \cP_{\Gr}$. This category is
endowed with a functor $\Res_d: Rep(G)\to \cA_d$ defined by $\Res_d(\cF )=
L_d*Z(\cF ) \mod \cP^I_{<\bc}$. The functor $\Res_d$ has natural automorphism
$\fM_d$ induced by the automorphism $\fM$ of
monodromy. The following Theorem is proved in \cite{B} Theorems 1 and 2:

{\bf Theorem.} (a) {\em The category $\cA_d$ has a natural structure of
a tensor category with unit object $L_d$, functor $\Res_d$ has a natural
structure of a tensor functor and $\fM_d$ is an automorphism of
the tensor functor $\Res_d$.}

(b) {\em Moreover, there exists a subgroup $H_d\subset G$, a unipotent
element $N_d\in G$ commuting with $H_d$, an equivalence of tensor
categories $\Phi_d: Rep(H_d)\to \cA_d$, and a natural transformation of
functors $\Res^G_{H_d}\simeq \Phi_d\circ \Res_d$ which interwines the tensor
automorphism $\fM_d$ with the action of $N_d$. The pair $(H_d, N_d)$ is unique
up to a simultaneous conjugacy. The element $N_d$ is conjugate to $u_{\bc}$.}

It is proved in \cite{LX} that the
 intersection $\bc \cap \hW^f$ consists of a unique
{\em canonical} left cell which we will denote $\Gamma_{\bc}$ (recall that
$\hW^f$ is a set of shortest representatives of right $W_f-$cosets in $\hW$).
In particular, there exists a
unique distinguished involution $d=d^f\in \bc \cap
\hW^f$. We call $d^f$ a {\em canonical} distinguished involution.

{\bf Theorem.} (see \cite{B} Theorem 3) (a) {\em The set of irreducible objects
of $\cA_{d^f}$ is $\{ L_w | w\in \Gamma_{\bc}\cap (\Gamma_{\bc})^{-1}\}$.}

(b) {\em The subgroup $H_{d^f}$ contains a maximal reductive subgroup of the
centralizer $Z_G(u_{\bc})$.}

\subsection{Central action of $Rep(F_{\bc})$}
Consider the functor $\tilde F: Rep(G)=\cP_{\Gr}\to \cA_{\bc}$ defined by
$\tilde F(\cF )=Z(\cF )*\fI_{\bc} \mod \cP^I_{<\bc}$. It is easy to see
from \ref{Gait} that the functor $\tilde F$ has a natural
 structure of a central functor.
Moreover, this functor has a canonical tensor unipotent automorphism
$\fM$ (monodromy) commuting with the centrality isomorphism.

{\bf Theorem 1.} {\em There exists a central functor $F: Rep(Z_G(u_{\bc}))\to
\cA_{\bc}$ such that $\tilde F=F\circ \Res^G_{Z_G(u_{\bc})}$. Moreover,
automorphism $\fM$ is induced by the action of $u_{\bc}$ on
$\Res^G_{Z_G(u_{\bc})}$.}

{\bf Proof.} Let $\cD (\cA_{\bc})$ denote the Drinfeld double of the monoidal
category $\cA_{\bc}$, see e. g. \cite{Kas}. By the universal property of double
the functor $\tilde F$ can be factorized as
$$
\xymatrix{Rep(G)\ar[rr]^{\tilde F} \ar[dr]^{\tilde F_0}&&\cA_{\bc}
\\ &\cD (\cA_{\bc})\ar[ur]&}
$$
where $\tilde F_0$ is again a central functor. Recall that the unit object of
$\cD (\cA_{\bc})$ is $\fI_{\bc}$ endowed with the centrality isomorphism
induced by the unity isomorphisms:
$$
\fI_{\bc}\bullet X\cong X\cong X\bullet \fI_{\bc}.
$$
We remark that $\fI_{\bc}$ considered as an object of $\cD (\cA_{\bc})$ is
irreducible: it is easy to see from Lemma \ref{peresechenie} that any
subobject of $\fI_{\bc}$ in $\cA_{\bc}$ does not lie in the center of
$\cA_{\bc}$ even on the level of $K-$theory. Now consider the
full subcategory
$\tilde \cD \subset \cD (\cA_{\bc})$ consisting of all subquotients of objects
$\tilde F_0(X),\; X\in Rep(G)$. Then the conditions of Proposition 1 \cite{B}
are satisfied for the pair $(\tilde \cD ,\tilde F_0)$. Consequently, the
functor $\tilde F_0$ factors through the restriction functor $\Res^G_H$
$$
\xymatrix{Rep(G)\ar[rr]^{\tilde F_0} \ar[dr]^{\Res^G_H}&&\tilde \cD
\\ &Rep(H)\ar[ur]^{F_0}&}
$$
for some subgroup $H\subset G$ and the action of
$\fM$ is given by some unipotent element $u\in Z_G(H)$. Theorem 2 of
\cite{B} identifies $u$ with $u_{\bc}$. Hence the subgroup $H$ is contained in
$Z_G(u_{\bc})$. Without loss of generality we can assume that
$H=Z_G(u_{\bc})$. We set the functor $F$ to be equal to the composition
$$
Rep(Z_G(u_{\bc}))\stackrel{F_0}{\longrightarrow} \tilde \cD \to \cD
(\cA_{\bc})\to \cA_{\bc}.
$$
The Theorem is proved. \sq

Let us restrict  $F$ to the semisimple part of the category
$Rep(Z_G(u_{\bc}))$, that is to the category $Rep(F_{\bc})$ where $F_{\bc}$ is
the maximal reductive factor of $Z_G(u_{\bc})$.

{\bf Proposition.} {\em For any $X\in Rep(F_{\bc})$ the object $F(X)\in
\cA_{\bc}$ is semisimple.}

{\bf Proof.} We can assume that $X$ is simple. Let $Y\in Rep(G)$ be an object
such that $X$ is a subquotient of $\Res^G_{F_{\bc}}(Y)$. The object
$\tilde F(Y)$ carries the monodromy filtration; by  Gabber's Theorem
(see \cite{BB}, Theorem 5.1.2) it coincides with the weight filtration, so
the  associated
graded object $gr \tilde F(Y)$ is semisimple by \cite{BBD}, 5.4.6.
 By the
Theorem 1 we get the same filtration from the action of $u_{\bc}$ on
$F(\Res^G_{Z_G(u_{\bc})})$. But the object $X$ is a direct summand of
$gr\; \Res^G_{Z_G(u_{\bc})}(Y)$ with respect to this filtration. \sq

As a corollary we get

{\bf Theorem 2.} {\em The functor $F$ restricts to a central functor
$F: Rep(F_{\bc})\to \cM_{\bc}$.}

\section{Canonical cell}

\subsection{Module categories} \label{mc}
In this subsection we review basic theory of module categories. A more
detailed exposition will appear in \cite{O}. We will work over a fixed field
$k$.

Let $\cC$ be an abelian monoidal category with biexact tensor product and with
 unit object $\be$.

{\bf Definition.} {\em A module category $\cM$ over $\cC$ is an abelian
category $\cM$ endowed with}

1) {\em An exact bifunctor $\otimes : Rep(H)\times \cC \to \cC$},

2) {\em Functorial} associativity isomorphisms {\em
$V\otimes (V'\otimes M)\simeq (V\otimes V')\otimes M$
for any $V, V'\in \cC, M\in \cM$},

3) {\em Functorial} unit isomorphisms {\em $\be \otimes M\to M$ for any $M\in
\cM$}

{\em subject to the usual pentagon and triangle axioms: the following
diagrams where all arrows are associativity and unit isomorphisms commute:}

$$
\xymatrix{&((V_1\otimes V_2)\otimes V_3)\otimes M \ar[dl] \ar[dr]&\\
(V_1\otimes (V_2\otimes V_3))\otimes M \ar[d]&&(V_1\otimes V_2)\otimes
(V_3\otimes M) \ar[d]\\V_1\otimes ((V_2\otimes V_3)\otimes M) \ar[rr]&&
V_1\otimes (V_2\otimes (V_3\otimes M))}
$$

$$
\xymatrix{(V\otimes \be )\otimes M\ar[rr] \ar[dr]&&V\otimes (\be \otimes M)
\ar[dl]\\ &V\otimes M&}
$$
The notions of module functors, and, in particular, equivalences of module
categories is defined in the obvious way.

{\bf Remark.} Module categories over general monoidal categories were
considered by L.~Crane and I.~Frenkel, see \cite{FC}. The name comes from
considering the notion of a
 monoidal category as categorification of the notion of a ring.
Module categories seem to be of importance in Conformal Field Theory where they
are implicitly considered in the context of Boundary Conformal Field Theory.

Of course the category $\cC$ is a module category over itself with
associativity and unit isomorphisms induced by ones in tensor category $\cC$.
Another example can be obtained as follows. Let $A\in \cC$ be an associative
algebra with unit, that is associative multiplication $A\otimes A\to A$ is
defined and there is an inclusion $\be \to A$ satisfying unit axioms. Then
category $Mod_{\cC}(A)$ of {\em right} $A-$modules in the category $\cC$ has
an obvious structure of a module category.

We will say that a module category $\cM$ is generated by objects $M_1, M_2,
\ldots \in \cM$ over $\cC$ if any object of $\cM$ is a subquotient of
$V\otimes M_i$ for some $V\in \cC$. We will say that $\cM$ is finitely
generated over $\cC$ if there exists finitely many (equivalently one)
objects $M_1, \ldots \in \cM$ such that $\cM$ is generated by them over
$\cC$.

Assume from now on that the category $\cC$ is rigid.
 Then there exists a canonical
isomorphism $Hom(V\otimes M, N)\cong Hom(M, V^*\otimes N)$ for any $V\in \cC ,
M, N\in \cM$.

Now assume that both categories $\cC$ and $\cM$ are semisimple. For any two
objects $M, N\in \cM$ the functor $\cC \to Vec_k,\; V\mapsto Hom(V\otimes M,N)$
is representable by 
an ind-object $\iHom (M, N)$ of $\cC$. By Yoneda's Lemma
$\iHom$ is a bifunctor $\cM^{op}\times \cM \to$ ind-objects of $\cC$.

{\bf Lemma.} {\em Assume that $\omega :\cM \to \cC$ is an exact faithful tensor
functor. Then for any} $M, N\in \cM\; \iHom (M, N)\in \cC$.

{\bf Proof.} It is clear that the
map $\iHom (M, N)\to \iHom (\omega (M), \omega
(N))$ is an imbedding. \sq

Assume that for any $M, N\in \cM$ the ind-object $\iHom (M, N)$ is
an object of $\cC$. For any three objects $M, N, K\in \cM$ a functorial and
associative multiplication $\iHom (N, K)\otimes \iHom (M, N)\to \iHom (M, K)$
is defined (note that the order of factors is opposite to the intuitive one).
 In
particular, for any object $M\in \cC$ the object $\iHom (M, M)$ has a natural
structure of an associative algebra in 
$\cC$. Assume that $\iHom (M, X)
\ne 0$ for any $X \in \cM$, that is the category $\cM$ is generated by $M$
over $\cC$. It is easy to see that
the functor $F_M: \cM \to Mod_{\cC}(A),\; F_M(X)=\iHom (M, X)$ is a tensor
functor. Moreover, we claim that this functor is an equivalence of categories.
The proof is straightforward: first one shows that the functor $F_M$ induces
an isomorphism on Hom's for objects of the form $V\otimes M,\; V\in \cC$, and
then one uses the fact that any object of $\cM$ admits a resolution by
objects of the form $V\otimes M$. Summarizing we get the following

{\bf Proposition.} {\em Let $\cC$ be a semisimple rigid monoidal category
and let $\cM$ be a semisimple module category over $\cC$. Assume that
there exists an exact faithful module functor $\omega :\cM \to \cC$. Then
the category $\cM$ is equivalent to the category $Mod_{\cC}(A)$ for some
associative algebra $A$. Moreover one can choose} $A=\iHom (M, M)$ {\em
for any object $M\in \cC$ generating $\cM$ over $\cC$.}

Let $\cM =Mod_{\cC}(A)$ be a module category. Consider the category $Fun(\cM ,
\cM )$ consisting of module functors $\cM \to \cM$. It is clear that
the category $Fun(\cM , \cM )$ is a
monoidal category with tensor product induced
by the composition of functors and identity functor as unit object. One shows
easily that the
monoidal category $Fun(\cM , \cM )$ is equivalent to the category
of $A-A$ bimodules in 
$\cC$ with the obvious monoidal structure.

\subsection{Module categories over $Rep(H)$}
In this subsection we specialize ourselves to the case when $\cC =Rep(H)$
for some reductive group $H$ over an algebraically closed field $k$ of
characteristic zero.

{\bf Examples.} (i) $Rep(H)$ with the
associativity and unit isomorphisms induced
from those in the monoidal category $Rep(H)$ is of course a
module category over $Rep(H)$.

(ii) More generally, let $X$ be a variety endowed with an
 $H-$action. The category
$Coh_H(X)$ of coherent $H-$equivariant sheaves on $X$ is a module category. We
get  example (i) by letting $X=$point.

(iii) Let $1\to \Gm\to \tilde H\to H\to 1$ be a central extension of $H$ whose
kernel is identified with
the multiplicative group $\Gm$
(we will call such a data just ``a central extension''); of course,
such an extension
is necessarily the pushforward if a central
extension $1\to C \to \tilde H'\to H\to 1$
under a homomorphism $C\to \Gm$ for a finite cyclic group $C$.
 Then the category $Rep^1(\tilde H)$
of representations $V$ of $\tilde H$ such that $\Gm$ acts on $V$ via identity
character is a module category over $Rep(H)$. We will also consider the
category
$Rep^{-1}(\tilde H)$ of representations of $\tilde H$ on which $\Gm$ acts
via character $x\mapsto x^{-1}$.

We will say that a module category $\cC$ has a {\em quasifiber functor} if
there exists a faithful exact module functor $\omega : \cC \to Rep(H)$.
The quasifiber functor if it exists is not unique: for any $V\in Rep(H)$ and
quasifiber functor $\omega$ the functor $M\mapsto \omega (M)\otimes V$ is
again a quasifiber functor.

{\bf Example.} (iv) Let $H'\subset H$ be a subgroup of finite index. Let
$1\to \Gm\to \tilde H'\to H'\to 1$ be a central extension.
 The category $Rep^1(\tilde H')$ is a module category over
$Rep(H)$ with the $Rep(H)-$action which factors through the restriction functor
$Rep(H)\to Rep(H')$. Let $V_0\in Rep^{-1}(\tilde H')$ be a fixed object. It is
easy to see that the functor $V\mapsto Ind_{H'}^H(V\otimes V_0), V\in
Rep^1(\tilde H')$ is a quasifiber functor ($\Gm$
 acts trivially on $V\otimes V_0$
so $V\otimes V_0$ can be considered as a representation of $H'$).

Example (iv) reduces to the example (i) with $X=H/H'$ if
the central extension splits.

{\bf Example.} (iv$'$) Finite sums of categories considered in Example (iv)
admit the following invariant description.

A finite {\it $H$-set of centrally extended points}
is the following collection of data:

(a) A finite set $X$ with an $H$ action;

(b) For any $x\in X$ a central extension $\Gm\to \Htil(x) \to H(x)$
of the stabilizer $H(x)=Stab_H(x)$.
These should be {\it equivariant} under the action of $H$, i.e. for
every $g\in G$ an isomorphism of $i_x^g:\Htil(x)\iso \Htil(gx)$ identical on
$\Gm$ and covering the map $C_g:H(x) \to H(g(x))$ (conjugation by $x$) should
be given. $i_g^x$ should coincide with the conjugation by $g$ when $g\in H(x)$
and should satisfy $i_{g_1g_2}^x=i_{g_1}^{g_2(x)}\circ i_{g_2}^x$.

Let $\bX$ be a finite set of centrally extended points. An {\it
equivariant sheaf}  on $\bX$ is a sheaf $\cF$
of finite dimensional $\BC$-vector
spaces on the underlying set $X$ together with

(a) a projective $H$-equivariant
structure on $\cF$.

(b) For every  $x\in X$ an action of $\Htil(x)$ on the stalk $\cF_x$,
comprising 
an object of $Rep^1(\Htil(x))$.

The data (a) and (b) should be compatible, i.e.
(b) should be $H$-equivariant, and the
 projective action of $H(x)$ arising from (b) must coincide
with the one arising from (a).

\medskip

Equivariant sheaves on $\bX$ obviously form a category, which we denote
by $Coh_H(\bX)$.

Choosing a set of representatives $x_i$ for $H$-orbits on $X$ we see
that the data of a centrally extended set with underlying equivariant
set $X$ is equivalent to a collection $\tilde H_i$ of central extensions
$\Gm\to \Htil_i\to H_i=H(x_i)$.  The category
$Coh_H(\bX)$ is then canonically equivalent to the direct sum $\oplus_i
Rep^1(\Htil_i)$.

{\bf Theorem 3.} {\em Let $\cM$ be a semisimple module category over $Rep(H)$
finitely generated over $Rep(H)$. Assume that $\cM$ admits
a quasifiber functor.
Then $\cM$ is equivalent to $Coh_H(\bX)$ for some centrally extended
finite $H$-set $\bX$
(i.e. to a
 finite direct sum of some categories of the type described in Example} (iv)
 {\em above).}

{\bf Proof.} By Proposition \ref{mc} the module category $\cM$ is equivalent to
the module category $Mod_{Rep(H)}(A)$ for some finite dimensional $H-$algebra
$A$.

{\bf Lemma.} {\em Semisimplicity of $\cM$ implies semisimplicity of $A$ as
an algebra in the category of vector spaces.}

{\bf Proof.} Consider the
regular representation $A_{reg}$ of $A$ as an object of
$Mod_{Rep(H)}(A)$. Let $r(A)$ be the Jacobson radical of $A$. It is clear that
$r(A)$ is $H-$invariant, hence $r(A)$ is subobject of $A_{reg}$ in
$Mod_{Rep(H)}(A)$.
Suppose $A_{reg}=r(A)\oplus A_1$ for some $A_1\in Mod_{Rep(H)}(A)$.
Applying the forgetful functor $Mod_{Rep(H)}(A)\to Mod(A)$ to $A_{reg}$ we
would get a complement to $r(A)$, which is impossible unless $r(A)=0$. \sq

Now let $A\ni 1=\sum e_i$ is the decomposition of $1$ in the sum of minimal
central orthogonal idempotents. The group $H$ acts on the set $\{ e_i\}$. We
may assume that this action is transitive. Let $H_1\subset H$ be the stabilizer
of $e_1$, the subgroup of finite index in $H$. The algebra $e_1Ae_1$ is
isomorphic
to the matrix algebra and the group $H_1$ acts on $e_1Ae_1$. We can choose
a projective representation $V$ of $H_1$ and an
isomorphism $e_1Ae_1\cong End(V)$.
It is clear that $A\cong Ind_{H_1}^H(e_1Ae_1)=Ind_{H_1}^HEnd(V)$.

The projective action of $H_1$ on $V$ comes from an action of a central
extension $\tilde H_1$ of $H_1$. Let us consider the corresponding category
from example (iv) $Rep^1(\tilde H_1)$. The representation $V$ can be viewed as
an object of this category and one easily calculates $\iHom (V, V)=
Ind_{H_1}^HEnd(V)$. The Theorem is proved. \sq


\subsection{Module category corresponding to the canonical cell}
For any subset $S\subset \bc$ let $\cM_S\subset \cM_{\bc}$ denote the full
Serre subcategory with simple objects $L_w, \; w\in S$.
Let $\Gamma \subset \uc$ be the canonical right cell, see \cite{LX}.
Let $\cM_{\Gamma}\subset \cM_{\uc}$ be the corresponding subcategory.
By the definition of a right cell we have $\cM_{\Gamma}\bullet \cM_{\bc}\subset
\cM_{\Gamma}$. Define on $\cM_{\Gamma}$ a structure of a module category
over $Rep(F_{\bc})$ by the formula $V\otimes M=F(V)\bullet M$ where $F$ is
a functor from Theorem 2. Note that due to the centrality of functor $F$
we have $F(Rep(F_{\bc}))\bullet \cM_{\Gamma}=\cM_\Gamma \bullet F(Rep(F_{\bc}))
\subset \cM_\Gamma$ so this is well defined. We claim that this category
admits a quasifiber functor. Indeed, let $\{ w_1, w_2, \ldots \} \subset
\Gamma^{-1}$ be a set of representatives of all {\em right} cells contained in
$\uc$ (such a set exists by the Lemma \ref{peresechenie} and is finite since
Lusztig proved (see \cite{Lc} II 2.2) that the number of cells in an affine
Weyl group is finite). Consider the functor $\cM_{\Gamma}\to \cM_{\Gamma \cap
\Gamma^{-1}},\; X\mapsto X\bullet (\oplus L_{w_i})$. Recall that in \cite{B}
the monoidal category $\cM_{\Gamma \cap \Gamma^{-1}}$ was identified with
$Rep(F_{\bc})$, see \ref{Bezr}. It is a simple exercise to check that this
functor is module functor with the module structure induced by the
associativity isomorphism in $\cM_{\bc}$, and it is clear that it is exact and
faithful. So this is quasifiber functor, and we can apply Theorem 3. We get

{\bf Proposition.} {\em The category $\cM_{\Gamma}$ as a module category over
$Rep(F_{\bc})$ is equivalent to the category $Coh_{F_{\bc}}(\bX )$ of coherent
sheaves on a finite $F_{\bc}-$set $\bX$ of possibly centrally extended points.}

Note that the inclusion $\cM_{\Gamma \cap \Gamma^{-1}}\subset \cM_\Gamma$ gives
us a distinguished point $\bo \in \bX$ which is just a usual (not centrally
extended) point fixed by the $F_{\bc}-$action.

\section{Square of a finite set}
\subsection{Monoidal category $Fun_{F_{\bc}}(\bX ,\bX )$}
Consider the category $Fun_{F_{\bc}}(\bX , \bX )$ consisting of all module
functors $Coh_{F_{\bc}}(\bX )\to Coh_{F_{\bc}}(\bX )$. It is a monoidal
category with the tensor product induced by the composition of functors and
 unit object equal to the identity functor. Since the category $Coh_{F_{\bc}}
(\bX)$ is semisimple, any functor $F\in Fun_{F_{\bc}}(\bX , \bX )$ has left
and right adjoint functors $F^*$ and $^*F$. Observe that adjoint of tensor
functor has a natural structure of a module functor and hence $F^*, {}^*F\in
Fun_{F_{\bc}}(\bX , \bX )$. Standard properties of adjoint functors show that
$F^*$ and $^*F$ are right and left duals of $F$ in the monoidal category
$Fun_{F_{\bc}}(\bX , \bX )$. So the category $Fun_{F_{\bc}}(\bX , \bX )$ is
rigid.

{\bf Lemma.} {\em The category $Fun_{F_{\bc}}(\bX , \bX )$ is semisimple.}

{\bf Proof.} Let us choose an $F_{\bc}-$algebra $A$ and an equivalence
$Coh_{F_{\bc}}(\bX )\to Mod_{Rep(F_{\bc}}(A)$. Then the category $Fun_{F_{\bc}}
(\bX , \bX )$ is equivalent to the category of $A-A$ bimodules in
$Rep(F_{\bc})$, or to the category of $A\otimes A^{op}-$modules in
$Rep(F_{\bc})$ where $A^{op}$ is $A$ with the opposite multiplication.
The latter category is clearly semisimple since $A\otimes A^{op}$ is
a semisimple algebra. \sq

Note that in semisimple monoidal category left and right duals coincide, so
in the future we will not distinguish left and right duals.

{\bf Remark.} For an $H$-set $X$ it is easy to construct an
 equivalence $Fun_{H}(X , X )\cong Coh_H(X\times X)$. Let us spell
out a generalization of this statement to centrally extended $H$-sets.

Recall that
for two central extensions $1\to \Gm\to \Htil_i \to H\to 1$, $i=1,2$
their product is defined by $\Htil_{12}=\Htil_1\times_H \Htil_2
/\Gm$, where $\Gm$ is embedded antidiagonally; also for  a central
extension $\Htil$ the
opposite central extension $\Htil'$
is the same group with the same homomorphism to
 $H$ but with the
identification of its kernel with $\Gm$ replaced by the opposite one
(composition of the original one with
 the map $x\mapsto x^{-1}$).

Now for two centrally extended $H$-sets $\bX$, $\bY$
one can define their product in the obvious manner:
the  underlying equivariant set is $X\times Y$, where $X,Y$
are equivariant sets  underlying
 $\bX$ and $\bY$; the central extension $\Htil(x,y)$ is the product
of restrictions of $\Htil(x)$ and of $\Htil(y)$ to $H(x,y)$.
 For a centrally extended $H$-set $\bX$ we obtain the {\it opposite}
centrally extended set $\bX'$
replacing each of the central extensions $\Gm\to \Htil(x)\to H(x)$
by the opposite one.

If $\bX$, $\bY$, $\bZ$ are centrally extended $H$-sets
with underlying $H$-sets $X,Y,Z$, then
for $\cF\in Coh_H(X\times Y')$, $\cG\in Coh_H(Y\times Z)$ the
sheaf $\cF\boxtimes \underline\BC\otimes \underline\BC\otimes \cG$
on $X\times Y\times Z$ carries a natural structure of an equivariant
sheaf on $\bX\times \bY_{0}\times \bZ$ (here $\bY_0$ is $Y$
equipped with the trivial (split) centrally
extended structure). Thus  the {\it convolution}
$\cF*\cG=\pr_{13*}(\cF\boxtimes \underline\BC\bigotimes
 \underline\BC\boxtimes \cG)$
(where $\pr_{13}:X\times Y\times Z\to X\times Z$ is the projection)
carries the structure of an equivariant sheaf on  $\bX\times \bZ$.
In particular, for $\bX=\bY=\bZ$ we get a monoidal structure on $Coh
_H(\bX\times \bX')$; and for $\bX=\bY$, and $\bZ$ being the point with
the split central extension we get a monoidal functor of
$Coh_H(\bX\times \bX')\to Fun_H(\bX)$. It is easy to see that this
functor is an equivalence.


\subsection{Monoidal functor $G$} \label{main}
We have a monoidal functor
$G: \cM_{\bc}^{op}\to Fun_{F_{\bc}}(\bX ,\bX ),$
$G(X)=?\bullet X$ where
$\cM_{\bc}^{op}$ is $\cM_{\bc}$ with the
{\em opposite} tensor product.
It is clear that  $G$ is exact and faithful.

The main result of this section is the following

{\bf Theorem 4.} {\em The functor $G$ is a
tensor equivalence $\cM_{\bc}^{op}\to Fun_{F_{\bc}}(\bX , \bX )$.}

{\bf Corollary.} {\em Suppose that
any subgroup of finite index in $F_{\bc}$ has no nontrivial
projective representations. Then Lusztig's Conjecture holds for the cell
$\uc$.}

\subsection{A result of G.~Lusztig} \label{dist}
The following result cited from  \cite{Lc} II Proposition 1.4
is cruicial for the proof of Theorem 4.

{\bf Proposition.} (a) {\em Assume that $L_x\bullet L_y,\; x, y\in \bc$
contains as a direct summand $L_d,\; d\in \cD$. Then $x=y^{-1}$ and
the multiplicity of $L_d$ in $L_x\bullet L_y$ is one.}

(b) {\em For any $x\in \bc$ the truncated convolution $L_x\bullet L_{x^{-1}}$
contains $L_d$
for a uniquely defined $d\in \cD \cap \bc$.}

\subsection{Proof of the Theorem 4.}
Since  the category
$Fun_{F_{\bc}}(\bX , \bX )$ is semisimple it is enough to
prove the following statements:

(i) Any functor from $Fun_{F_{\bc}}(\bX , \bX )$ appears as a direct summand
of $G(L),\; L\in \cM_{\bc}^{op}$.

(ii) For any $w\in \bc$ the functor $G(L_w)$ is irreducible.

(iii) For $w, w'\in \bc$ an isomorphism $G(L_w)=G(L_{w'})$ implies $w=w'$.

\subsubsection{} \label{step1}
We begin with the following

{\bf Lemma.} (a) {\em For any $w\in \Gamma$ the functor $G(L_w)$ is
irreducible. Moreover $G$ induces an equivalence $\cM_{\Gamma}\to$
\{ module functors $Rep(F_{\bc})\to \cM_{\Gamma}$\} .}

(b) {\em For any $w\in \Gamma^{-1}$ the functor $G(L_w)$ is irreducible.
Moreover, $G$ induces an equivalence $\cM_{\Gamma^{-1}}\to$ \{ module functors
$\cM_{\Gamma}\to Rep(F_{\bc})$\} .}

{\bf Proof.} (a) It is clear that $L_v\bullet L_w=0$ for any $w\in \Gamma ,\;
v\in \Gamma -\Gamma \cap \Gamma^{-1}$. So the functor $G(L_w)$ can be
considered as a functor $Rep(F_{\bc})=\cM_{\Gamma \cap \Gamma^{-1}}\to
\cM_{\Gamma}$. For any module category $\cM$ over $Rep(F_{\bc})$ the map
$f\mapsto f(\be )$ defines an equivalence of categories $Fun_{Rep(F_{\bc})}
(Rep(F_{\bc}), \cM)\to \cM$. In our case
$G(L_w)(\be )=L_{d^f}\bullet L_w=L_w$ and (a) is proved.

(b) Let us first check that for $w\in \Gamma^{-1}$ we have
$$
G(L_w)^*\cong G(L_{w^{-1}}),
\eqno{(*)}
$$
where $G(L_w)^*$ is the functor adjoint to $G(L_w)$.
 
For $w\in \Gamma^{-1}$ the functor $G(L_w)$ maps
$\cM_{\Gamma}$ to $\cM_{\Gamma \cap \Gamma^{-1}}=Rep(F_{\bc})\subset
\cM_\Gamma$. Thus $G(L_w)^*$ sends $L_v$ to zero
unless $v\in \Gamma\cap \Gamma^{-1}$; hence part (a) of the Lemma implies
that $G(L_w)^*$ is isomorphic to $G(L)$ for some $L\in \cM_\Gamma$.
To see that 
$
L\cong L_{w^{-1}}
$
 it is enough to check that for $v\in \Gamma$ the space
$Hom((G(L_v),G(L_w)^*))$
 is one dimensional if $v=w^{-1}$, and is zero otherwise.
We have $Hom((G(L_v),G(L_w)^*))=Hom(G(L_w)\circ G(L_v),Id_{\cM_\Gamma})$.
Notice that $ G(L_w)\circ G(L_v)$ preserves the direct summand 
$\cM_{\Gamma \cap \Gamma^{-1}}\subset \cM_\Gamma$ and is zero on its 
complement. It follows that $$Hom(G(L_w)\circ G(L_v),Id_{\cM_\Gamma})=
Hom(G(L_w)\circ G(L_v)|_{\cM_\Gamma\cap \Gamma^{-1}}
,Id_{\cM_\Gamma\cap \Gamma^{-1}}).$$
Since $\cM_{\Gamma\cap \Gamma^{-1}}\cong Rep(F_{\bc})$ by \ref{Bezr},
and  the category of module functors from $Rep(H)$ to $Rep(H)$
(considered as the free module over itself) is equivalent to $Rep(H)$,
we see that $$Hom (G(L_w)\circ G(L_v),Id_{\cM_{ \Gamma\cap \Gamma^{-1}} })
= Hom_{\cM_{\Gamma\cap \Gamma^{-1}} } (L_w\bullet L_v, L_{d^f}),$$
thus $(*)$  follows from Proposition \ref{dist}.

Irreducibility of $G(L_w)$ follows from $(*)$  and part (a), because
 the dual object of an irreducible object (in the category of functors)
is irreducible.
It remains to check that any module functor $\phi:\cM_\Gamma
\to Rep(F_{\bc})$ is isomorphic to the one coming from some 
$L\in \cM_{\Gamma^{-1}}$. Consider $\phi$ as an endofunctor $\cM_\Gamma$
(i.e. take its composition with the imbedding $Rep(F_{\bc})=\cM_{\Gamma
\cap \Gamma^{-1}}\imbed \cM_\Gamma$); then by (a) the adjoint
functor $\phi^*$ is isomorphic to $G(L)$ for $L\in \cM_\Gamma$,
so this statement also follows from $(*)$.
\sq

We can now prove (i).

{\bf Corollary.} {\em Any irreducible functor from $Fun_{Rep(F_{\bc})}(\bX ,
\bX )$ appears as a direct summand of $G(L_{w'}\bullet L_w), w\in \Gamma ,
w'\in \Gamma^{-1}$.}

{\bf Proof.} We need to prove that any irreducible functor $f\in
Fun_{Rep(F_{\bc})}(\bX , \bX )$ is a direct summand of a composite functor
$Coh_{F_{\bc}}(\bX )\to Rep(F_{\bc})\to Coh_{F_{\bc}}(\bX )$. For this
we choose any functor $g : Coh_{F_{\bc}}(\bX )\to Rep(F_{\bc})$ such that
the composition $g\circ f$ is nonzero. Let $g^*$ be the adjoint functor to $g$.
Then $f$ is evidently a direct summand of $(g^*\circ g) \circ f$ which
admits the required factorisation $g^*\circ (g\circ f)$. \sq

\subsubsection{} \label{step2}
{\bf Lemma.} {\em For any $w\in \uc -\cD$ we have $Hom(G(L_w), id)=0$.}

{\bf Proof.} We note that $L_w$ as an object of $\cM_{\bc}$ is a direct summand
of $L_u\bullet L_v$ where $u\in \Gamma^{-1}$ and $v\in \Gamma$ (this follows
easily for the example from
 Theorem 1.8 in \cite{Lc} II). Assume that
$Hom(G(L_w), id)\ne 0$ and hence $Hom(G(L_u\bullet L_v), id)\ne 0$.
Then we get a nonzero transformation $G(L_u)\to G(L_v)^*$;
since both functors are irreducible by  Lemma \ref{step1} they are actually
isomorphic. On the other hand,  Proposition \ref{dist} provides
a non-zero transformation $G(L_{u})\circ G(L_u^{-1})\to Id$, which also
yields an isomorphism $G(L_u)^*\cong G(L_{u^{-1}})$. Thus
$G(L_{u^{-1}})\cong G(L_v)$, and by Lemma
\ref{step1} this yields $L_{u^{-1}}\cong L_v$, so $v=u^{-1}$.
 Furthermore,
$dim Hom(G(L_u\bullet L_v), id)=dim Hom (G(L_v),
G(L_u)^*)=1$; and  by
Proposition \ref{dist} the object $L_u\bullet L_v$ contains $L_d$ for a
 uniquely
defined $d\in \cD$. Since $Hom(G(L_d), id)\ne 0$ we have that
$Hom(L_w, id)=0$ if $w\ne d$. The Lemma is proved. \sq

Now we can prove (ii)

{\bf Corollary.} {\em For any $w\in \uc$ the functor $G(L_w)$ is irreducible.}

{\bf Proof.} Consider the adjoint functor $G(L_w)^*$. By  Lemma \ref{step1}
any summand of $G(L_w)^*$ appears as a direct summand of $G(L_{w'})$.
For such $w'$ we have a non-zero transformation $G(L_{w'})\circ
G(L_w)=G(L_w \bullet L_{w'})\to Id$, and by Lemma \ref{step2} and Proposition
\ref{dist} this is possible only when $w'=w^{-1}$. If $G(L_w)$ is reducible
this implies that $dim Hom(G(L_w)\circ G(L_{w^{-1}}), id)>1$. On the other hand
$dim Hom(G(L_w)\circ G(L_{w^{-1}}), id)=dim Hom(G(L_{w^{-1}}
\bullet L_w), id)=1$ by Proposition \ref{dist} and  Lemma \ref{step2},
and we get a contradiction. \sq

\subsubsection{} We can now prove (iii). Assume that $G(L_w)=G(L_{w'})$. Then
$G(L_w)\circ G(L_{w^{-1}})=G(L_{w'})\circ G(L_{w^{-1}})$. Since
$Hom(G(L_w)\circ G(L_{w^{-1}}), id)\ne 0$ we have that
$Hom(G(L_{w'})\circ G(L_{w^{-1}}), id)\ne 0$ and by the Proposition \ref{dist}
$w'=(w^{-1})^{-1}=w$.

Since we proved (i) (ii) and (iii) the Theorem 4 is proved. \sq

\subsection{Examples} The Corollary \ref{main} can be applied in the
following cases:

(a) Let $G$ be a simply connected group and let $\bc$ be the lowest cell.
In this case $u_{\bc}=e\in G$ and $F_{\bc}=G$. In this case Corollary
\ref{main} is a result of \cite{X1}.

(b) Let $G=GL_n$. In this case all groups $F_{\bc}$ are connected and have
no nontrivial projective representations (these groups are products of various
$GL_m$). In this case we get a result of \cite{X3}.

(c) Let $G$ be a simply connected group of rank 2. In this case one easily
verifies that the condition of Corollary \ref{main} is satisfied and we get a
result of \cite{X2}.

(d) Let $G$ be a simple simply connected group. Let $\bc$ be the subregular
cell, that is the
cell corresponding to the subregular nilpotent orbit. Again one
easily verifies that the condition of Corollary \ref{main} is satisfied except
if $G$ is of type $C_n$. In the latter case $F_{\bc}=\BZ /2\BZ \times
\BZ /2\BZ$ where one of the factors comes from the center of $G$. One can
exlude centrally extended points in this case by considering a
reductive group
$G_1=G\times T/(z, -1)$ where $T$ is the
one dimensional torus, $z\in G$ is the unique
nontrivial central element and $-1\in T$ is the unique nontrivial involutive
element. So we get another result of \cite{X2}.

Finally note that centrally extended points naturally appear in the 
description of truncated convolution categories for simple non 
simply-connected groups, see in \cite{X3} 8.3 example with $G=PSL_2$.

\end{document}